\documentclass[11pt,leqno]{article}
\usepackage{amsmath}  
\usepackage{amsfonts}  
\usepackage[cp850]{inputenc} 
\begin{document}  
\baselineskip=1.1\baselineskip 
\newcommand{\blockdiag}{\mbox{blockdiag}}  
\newcommand{\rank}{\mbox{rank}}  
\newcommand{\ds}{\displaystyle}  
\newcommand{\sss}{\scriptstyle}  
\newcommand{\R}{\mathbb{R}}  
\newcommand{\Naturals}{\mathbb{N}}  
\newcommand{\un}{\underline}  
\newcommand{\hs}{\hspace}  
\newcommand{\vs}{\vspace}  
\newcommand{\no}{\noindent}  
\newcommand{\be}{\begin{equation}}  
\newcommand{\ee}{\end{equation}}  
\newcommand{\ba}{\begin{array}}  
\newcommand{\ena}{\end{array}}  
\newcommand{\diag}{\mbox{\rm diag}}  
\newcommand{\1}{\mbox{\bf 1}}  
\newcommand{\bdes}{\begin{description}}  
\newcommand{\edes}{\end{description}}  
\newcommand{\bla}{\begin{lemma}}  
\newcommand{\ela}{\end{lemma}}  
\newcommand{\bcy}{\begin{corollary}}  
\newcommand{\ecy}{\end{corollary}}  
\newtheorem{thm}{Theorem}[section]   
\newtheorem{re}[thm]{Remark}  
\newtheorem{co}[thm]{Corollary}  
\newtheorem{pr}[thm]{Proposition}  
\newtheorem{de}[thm]{Definition}  
\newtheorem{lm}[thm]{Lemma}  
\newtheorem{example}[thm]{Example}  
\renewenvironment{description}
  {\begin{list}{}%
   {\itemsep=0pt%
   \itemindent=0pt%
   \listparindent=0pt%
   \labelwidth=20pt%
   \leftmargin=30pt}}%
  {\end{list}}%
\renewcommand{\theequation}{\thesection.\arabic{equation}}  
\def\squarebox#1{\hbox to #1{\hfill\vbox to #1{\vfill}}}  
\newcommand{\qed}{\hspace*{\fill}  
\vbox{\hrule\hbox{\vrule\squarebox{.45em}\vrule}\hrule}\smallskip}  
\def \ci {\mathop{\hbox {\vrule height .4pt depth 0pt width 0.5cm  
\hskip -0.3cm \vrule height 10pt  
depth 0pt \hskip 0.1cm \vrule height 10pt depth 0pt \hskip 0.2cm}}}  
\def\E{\mbox{\rm E}}  
\def\supp{\mbox{\rm supp}}  
\def\Cov{\mbox{\rm Cov}}  
\def\realpart{\mbox{\textsc{Re}}}  
%
\title{Linear stochastic differential-algebraic equations  
\\  
with constant coefficients  
\vs{1cm}}  
\author{  
Aureli Alabert  
\thanks{Supported by grants 2005SGR01043 of CIRIT and BFM2003-0261 of MCYT}  
\\  
Departament de Matem\`atiques \\  
Universitat Aut\`onoma de Barcelona \\  
08193 Bellaterra, Catalonia \\  
e-mail: alabert@mat.uab.cat  
\and  
Marco Ferrante  
\thanks{Partially supported by a grant of the CRM, Bellaterra, Spain
and a grant of the GNAMPA, Italy}  
\\  
Dipartimento di Matematica P. e A. \\  
Universit\`a degli Studi di Padova \\  
via Belzoni 7,  
35131 Padova, Italy \\  
e-mail: ferrante@math.unipd.it}  
\maketitle  
\begin{abstract}  
\no  
We consider linear stochastic differential-algebraic equations  
with constant coefficients and additive white noise.
Due to the nature of this 
class of equations, the solution must be defined as a generalised process
(in the sense of Dawson and Fernique). We provide sufficient
conditions for the law of the
variables of the solution process to be absolutely continuous with
respect to Lebesgue measure.
\end{abstract}  
  
\vspace{3truecm}  
  
\thispagestyle{empty}  
\vfill\eject  
\null  
\thispagestyle{empty}  
\vfill\eject  
\setcounter{page}{1}  
\section{Introduction}  

  A Differential-Algebraic Equation (DAE) is, essentially, an Ordinary 
  Differential Equation (ODE)
  $F(x,\dot x)=0$ that cannot be solved for the derivative $\dot x$. 
  The name comes from the fact that in some cases they can be 
  reduced to a two-part 
  system: A usual differential system plus a ``nondifferential'' one 
  (hence ``algebraic'', with some abuse of language), that is 
\begin{equation}
\label{1.1}
\left\{
\begin{aligned}
         \dot x_1&=f(x_1,x_2)
         \\
         0 &=g(x_1,x_2)
\end{aligned}
\right.
\end{equation}
  for some partitioning of the vector $x$ into vectors $x_1$ and $x_2$.
  In general, however, such a splitting need not exist.
\par 
  In comparison with ODE's, these equations present at least two
  major difficulties: the first lies in the fact that it is not possible to
  establish general existence and uniqueness results, due to their more
  complicate structure; the second one
  is that DAE's do not regularise the input (quite the contrary), 
  since solving them 
  typically involves differentiation in place of integration.
  At the same time, DAE's are very important objects, arising  
  in many application fields; 
  among them we mention 
  the simulation of electrical circuits, the modelling of 
  multibody mechanisms, the approximation of singular perturbation problems 
  arising e.g.~in fluid dynamics, the discretisation of partial differential 
  equations, the analysis of chemical processes, 
  and the problem of protein folding. 
  We refer to Rabier and Rheinboldt \cite{MR2003f:34002} for a survey
  of applications.
\par
  The class of DAE's most treated in the literature 
  is, not surprisingly, that of linear equations, which have the form 
\begin{equation*}
  A(t)\dot x(t)+B(t)x(t)=f(t)
  \ ,
\end{equation*}
  with $x,f\colon\R^+\rightarrow\R^n$ and 
  $A,B\colon \R^+\rightarrow\R^{n\times n}$. When $A$ and $B$ are constant
  matrices the equation is said to have constant coefficients. 
  Note that these equations cannot in general be split as in (\ref{1.1}).
\par
  Recently, there has been some incipient work (Schein and Denk \cite{MR99i:65076} and
  Winkler \cite{MR2004f:60133}) on Stochastic Differential-Algebraic 
  Equations (SDAE).
  In view to incorporate to the model a random external perturbation,
  an additional term is
  attached to the differential-algebraic equation, in the form of an 
  additive noise (white or coloured). 
  The solution will then be a stochastic process instead
  of a single function.

  Since the focus in \cite{MR99i:65076} and \cite{MR2004f:60133} 
  is on numerical solving and 
  the particular applications,
  some interesting theoretical questions have been left aside in these papers.
  Our long-term purpose is to put SDAE into the mainstream of stochastic calculus,  
  developing as far as possible a theory similar to that of 
  stochastic differential equations.
  In this first paper our aim is to
  investigate the solution of linear SDAE with constant coefficients  
  and an additive white noise, that means 
\begin{equation*}  
  A\dot x(t)+Bx(t)=f(t)+\Lambda \xi(t)
  \ ,
\end{equation*}
  where $\xi$ is a white noise and $A,B,\Lambda$ are constant matrices of 
  appropriate dimensions.
  We shall first reduce the equation to the so-called Kronecker Canonical Form (KCF),
    which is easy to analyse, and from whose solution one can recover
    immediately the solution to the original problem.
  Unfortunately, it is not possible to extend this approach to the case
  of linear SDAE with varying coefficients, just as happens in the deterministic case, 
  where several different 
  approaches have been proposed. Among these, the most
  promising in our opinion is that of Rabier and Rheinboldt \cite{MR97e:34011}.  
\par  
%
\par
  Due to the simple structure of the equations considered here,
  it is not a hard task to establish the existence of a unique solution 
  in the appropriate sense. However,
  as mentioned before, a DAE does not regularise the input $f(t)$ in general.
  If white noise, or a similarly irregular noise is used as input,
  then the solution process to a SDAE will not be a usual 
  stochastic process, defined as a random vector at every time $t$, but instead 
  a ``generalised process'', the random 
  analogue of a Schwartz generalised function. 
\par  
  The paper is organised as follows: in the next section we shall
  provide a short introduction to linear DAE's and to generalised processes.
  In the third section we shall define what we mean by a solution to a linear SDAE
  and in Section \ref{ac} we shall provide a sufficient condition for the
  existence of a density of the law of the solution. In the final Section \ref{ElecCirc} we
  shall discuss a simple example arising in the modelling of electrical circuits.
\par
  Superscripts in parentheses mean order of derivation. 
  The superscript $^\top$ stands for transposition.
  All function and vector norms
  throughout the paper will be $L^2$ norms, and the scalar product will be denoted
  by $\langle \cdot, \cdot \rangle$ in both cases.
  Covariance matrices of random vectors will be denoted by $\Cov(\cdot)$.
  The Kronecker delta notation $\delta_{ij}:=1_{\{i=j\}}$ will be used throughout.

\setcounter{equation}{0}  
\section{Preliminaries on DAE and generalised processes}  
\label{pre}  
  In this section we briefly introduce two topics: 
  the (deterministic) differential-algebraic equations   
  and the generalised processes. 
  An exhaustive introduction on the first topic can be found
  in Rabier and Rheinboldt \cite{MR2003f:34002}, while
  the basic theory of generalised processes 
  can be found in Dawson \cite{MR41:6332},
  Fernique \cite{MR0221576}, or Chapter 3 in Gel'fand and Vilenkin \cite{MR0435834}.

\subsection{Differential-Algebraic Equations}
\label{DAE}
  Consider an implicit autonomous ODE, 
\begin{equation}  
\label{ImpODE}
  F(x,\dot x)=0 \ , 
\end{equation}  
  where $F :=F(x,p):\R^{n\times n} \rightarrow \R^n$ is a 
  sufficiently smooth function. 
  If the partial differential $D_p F(x,p)$ is invertible at every point 
  $(x_0,p_0)$, one
  can easily prove that the implicit ODE is locally reducible to an
  explicit ODE. 
  If $D_p F(x_0,p_0)$ is not invertible, two cases are possible:
  either the total derivative $D F(x_0,p_0)$ is onto $\R^n$ or it is not.
  In the first case, and assuming that the rank of $D_p F(x,p)$ is
  constant in a neighbourhood of $(x_0,p_0)$, (\ref{ImpODE}) is called 
  a {\em differential-algebraic equation\/}, while in the remaining cases
  one speaks of an ODE with a singularity at $(x_0,p_0)$.

  A {\em linear DAE\/} is a system of the form
\begin{equation}  
\label{a-lin-SDAE}  
  A(t)\dot x + B(t) x=f(t)  
  \ ,\quad t\ge 0\ ,
\end{equation}  
  where $A(t), B(t)\in \R^{n \times n}$ and $f(t)\in \R^n$. 
  The matrix function $A(t)$ is assumed to have a constant (non-full) rank
  for any $t$ in the interval of interest. 
  (Clearly, if $A(t)$ has full rank for all $t$ in an interval, 
  then the DAE reduces locally to an ODE.)
  In the simplest case, when $A$ and $B$ do not depend on $t$, we have
  a {\em linear DAE with constant coefficients}, and 
  an extensive study of these problems has been developed.
  Since we want to allow solutions of DAE in the distributional sense, let us
  precise the definition of a solution.
  
  Let ${\cal D}'$ be the space of distributions
  (generalised functions) on some open set $U\subset\R$, that is, 
  the dual of the space ${\cal D}={\cal C}_{\text{c}}^{\infty}(U)$ of smooth functions with
  compact support defined on $U$. 
  An $n$-dimensional distribution is an element of $({\cal D}')^n$, and, for 
  $x=(x_1,\dots,x_n)\in({\cal D}')^n$ and $\phi\in {\cal D}$, we denote 
  by $\langle x,\phi\rangle=(\langle x_1,\phi\rangle,\dots,\langle x_n,\phi\rangle)^{\top}$
  the action of $x$ on $\phi$.
  We will always assume, without loss of generality, $U=]0,+\infty[$.
  
\begin{de}  
\label{solDAE}
  Let $f$ be an $n$-dimensional distribution on $U$, 
  and $A$, $B$ two $n\times n$
  constant matrices. A solution to the linear DAE with constant coefficients
\begin{equation}\label{linDAEcc}  
  A\dot x + Bx=f
\end{equation}  
  is an $n$-dimensional distribution $x$ on $U$ such that, for every test function 
  $\phi\in{\cal D}$, the following equality holds:
$$  
  A\langle \dot x,\phi\rangle 
  +
  B\langle x,\phi\rangle
  =
  \langle f,\phi\rangle
$$  
\end{de}
  The theory of linear DAE starts with the definition of a regular matrix pencil: 

  \begin{de}
\label{pencil}
  Given two matrices $A,B\in\R^{n\times n}$, 
  the {\em matrix pencil $(A,B)$} is the function 
  $\lambda\mapsto \lambda A+B$, for $\lambda\in\R$.
  It is called a {\em regular matrix pencil} if $\det(\lambda A+B)\neq 0$
  for some $\lambda$.
\end{de}
  If the matrices $A$ and $B$ in equation 
  (\ref{linDAEcc})
  form a regular matrix pencil, then a solution exists.
  This is a consequence of the following classical result due to 
  Weierstrass and Kronecker, which states that $A$ and $B$ can be
  simultaneously transformed into a convenient
  {\em canonical form}
  (see e.g. Griepentrog and M\"arz \cite{MR90d:34004} for the proof).
\begin{lm}
\label{kronecker}
  Given a regular matrix pencil $(A,B)$, there exist nonsingular $n\times n$ matrices
$P$ and $Q$ and integers $0\le d,q\le n$, with $d+q=n$, such that  
\begin{equation*}  
  PAQ=  
  \left(  
  \begin{matrix}  
  I_{d} & 0 \\  
  0 & N  
  \end{matrix}  
  \right)  
  \mbox{\quad and\quad}  
  PBQ=  
  \left(  
  \begin{matrix}  
  J & 0 \\  
  0 & I_{q}  
  \end{matrix}  
  \right)  
\end{equation*}  
  where $I_d, I_q$ are identities of dimensions $d$ and $q$,  
  $N=\blockdiag(N_1,\dots,N_r)$, with  
  $N_i$  
  the $q_i\times q_i$ matrix  
\begin{equation*}  
  N_i=  
  \left(  
  \begin{matrix}  
  0 & 1 & 0 & \dots & 0 \\  
  0 & 0 & 1 & \cdots & 0 \\  
  \vdots & & \ddots & \ddots & \vdots\\  
  \vdots & & & \ddots & 1 \\  
  0 & \cdots &\cdots &\cdots & 0  
  \end{matrix}  
  \right)  
  \ ,
\end{equation*}  
  and $J$ is in Jordan canonical form. 
\end{lm}
\begin{pr}
  If $(A,B)$ is a regular pencil, $x$ is a solution to the DAE (\ref{linDAEcc})
  if and only if $y=Q^{-1}x$ is a solution to   
\begin{equation}  
\label{reducedDAE}
  PAQ \dot y + PBQ y =Pf
  \ ,
\end{equation}
  where $P$ and $Q$ are the matrices of Proposition \ref{kronecker}. 
\end{pr}
\noindent
\emph{Proof:}
  The result is obvious, since (\ref{reducedDAE}) is obtained from 
  (\ref{linDAEcc}) multiplying from the left by the invertible matrix $P$. \qed

\par  
  System (\ref{reducedDAE}) is said to be in \emph{Kronecker Canonical Form} (KCF) 
  and splits into two parts. The first one is 
  a linear differential system of dimension $d$, 
  and the second one is an  
  ``algebraic system''
  of dimension $q$.
  Denoting by $u$ and $v$ the variables in the first and the second part
  respectively, and by $b$ and $c$ the related partitioning of the 
  vector distribution $Pf$ 
  we can write the two systems as follows:
\begin{equation}  
\label{3.3new}  
  \left(  
  \begin{matrix}  
  \dot u_1 \\  
  \vdots \\  
  \dot u_d  
  \end{matrix}  
  \right)  
  +  
  J  
  \left(  
  \begin{matrix}  
  u_1 \\  
  \vdots \\  
  u_d  
  \end{matrix}  
  \right)  
  =  
  \left(  
  \begin{matrix}  
  b_1 \\  
  \vdots \\  
  b_d  
  \end{matrix}  
  \right)  
  \ ,
\end{equation}  
\begin{equation}  
\label{3.4new}  
  N  
  \left(  
  \begin{matrix}  
  \dot v_{1} \\  
  \vdots \\  
  \dot v_{q}  
  \end{matrix}  
  \right)  
  +  
  \left(  
  \begin{matrix}  
  v_{1} \\  
  \vdots \\  
  v_{q}  
  \end{matrix}  
  \right)  
  =  
  \left(  
  \begin{matrix}  
  c_{1} \\  
  \vdots \\  
  c_{q}  
  \end{matrix}  
  \right)  
  \ .
\end{equation}  
  We refer to $u$ as the {\em differential variables\/} and to $v$ as
  the {\em algebraic variables.}
\par  
  The differential system has a unique solution once an initial condition,
  i.e. the value of the solution at some suitable test function $\phi_0$, 
  is given. The function must have a nonvanishing integral   
  (see Schwartz \cite{MR0209834}, V.6). It can be assumed without any 
  loss of generality that
  $\int_0^\infty \phi_0=1$.
  
  On the other hand, system (\ref{3.4new}) consists of a number of decoupled blocks, 
  which are
  easily and uniquely solved by backward substitution, without the need of any additional
  condition. For instance, for the first block,
\begin{equation}  
\label{AL1det}  
  N_1  
  \left(  
  \begin{matrix}  
  \dot v_{1} \\  
  \vdots \\  
  \dot v_{q_1}  
  \end{matrix}  
  \right)  
  +  
  \left(  
  \begin{matrix}  
  v_{1} \\  
  \vdots \\  
  v_{q_1}  
  \end{matrix}  
  \right)  
  =  
  \left(  
  \begin{matrix}  
  c_1 \\  
  \vdots \\  
  c_{q_1}  
  \end{matrix}  
  \right)
  \ ,
\end{equation}  
a recursive calculation gives the following distributional solution:
\begin{equation}\label{sae.dist.det}
  \langle v_j,\phi\rangle
  =
  \sum_{k=j}^{q_1}  
  \big\langle 
  c_k,\phi^{(k-j)}
  \big\rangle
  \ ,
  \quad j=1,\dots, q_1
  \ ,
  \quad \phi\in {\cal D}
\end{equation}  
  We can thus state the following proposition and corollary: 
\begin{pr}
\label{ExUnReduced}
  Assume $(A,B)$ is a regular matrix pencil. Then, for every $u^0=(u^0_1,\dots,u^0_d)\in\R^d$, 
  and every fixed test function
  $\phi_0$ with $\int_0^{\infty}\phi_0=1$, there exists a unique solution $y$ to (\ref{reducedDAE}) such that
\begin{equation*}  
  \langle y,\phi_0\rangle_i=u^0_i
  \ ,
  \quad
  i=1,\dots,d
  \ .
\end{equation*}
\end{pr}
  
\begin{co}
\label{ExUnOriginal}
  Assume $(A,B)$ is a regular matrix pencil. Then, for every $u^0=(u^0_1,\dots,u^0_d)\in\R^d$, 
  and every fixed test function
  $\phi_0$ with $\int_0^{\infty}\phi_0=1$, there exists a unique solution $x$ to (\ref{linDAEcc}) such that
\begin{equation*}  
  \langle Q^{-1}x,\phi_0\rangle_i=u^0_i
  \ ,
  \quad
  i=1,\dots,d
  \ .
\end{equation*}
\end{co}
\par
  Note that the matrix $N$ is nilpotent, with nilpotency index given by
  the dimension of its largest block. 
  The nilpotency index of $N$ in this
  canonical form is a characteristic of the matrix pencil and we shall call it the
  {\em index} of the equation (\ref{linDAEcc}). 
  The regularity of the solution depends directly on the index of the equation.

\begin{re}  
  Without the hypothesis of regularity of the pencil, a linear DAE may possess    
  an infinity of solutions or no solution at all, depending on the right-hand side.
  This is the case, for instance, of
\begin{equation*}  
  \left(  
  \begin{matrix}  
  0 & 0 \\  
  1 & 1  
  \end{matrix}  
  \right)
  \dot x(t)
  +
  \left(  
  \begin{matrix}  
  1 & 1 \\  
  0 & 0  
  \end{matrix}  
  \right)
  x(t)
  =
  \left(  
  \begin{matrix}  
  f_1(t) \\  
  f_2(t)   
  \end{matrix}  
  \right)
\end{equation*}
  with any fixed initial condition.
\end{re}
\subsection{Generalised processes}
\label{GenProc}
  As before, let ${\cal D}'$ be the space of distributions
  on an open set $U$.
  A {\em random distribution\/} on 
  $U$, defined in the probability space $(\Omega, {\cal F}, P)$, is a measurable mapping
  $X\colon(\Omega, {\cal F})\rightarrow ({\cal D}', {\cal B}({\cal D}'))$, where  
  ${\cal B}({\cal D}')$ denotes the Borel $\sigma$-field, relative to the weak-$\star$ topology
   (equivalently, the strong 
  dual topology, see Fernique \cite{MR0221576}). 
  Denoting by 
  $\langle X(\omega), \phi\rangle$ the action of the distribution $X(\omega)\in\cal D'$ 
  on the test function $\phi\in \cal D$,
  it holds that the mapping 
  $\omega\mapsto \langle X(\omega), \phi\rangle$
  is measurable from $(\Omega,\cal F)$ into $(\R,\cal B(\R))$, hence a real random 
  variable $\langle X,\phi\rangle$ on  $(\Omega,\cal F, P)$.
  The law of $X$ is determined by the law of the finite-dimensional vectors 
  $(\langle X,\phi_1\rangle,\dots,\langle X,\phi_n\rangle)$, $\phi_i\in\cal D$, $n\in\Naturals$. 
\par
  The sum of random distributions $X$ and $Y$ on $(\Omega,\cal F,P)$, defined in 
  the obvious manner, is again a random distribution. The product of a real random 
  variable $\alpha$ and a random distribution, defined by 
  $\langle \alpha X, \phi\rangle:=\alpha\langle X,\phi\rangle$, is also a random distribution.  
  The derivative of a random distribution, defined by 
  $\langle \dot X,\phi\rangle := -\langle X,\dot\phi\rangle$, is again a random distribution.
\par
  Given a random distribution $X$, the mapping $X\colon {\cal D}\rightarrow L^0(\Omega)$
  defined by $\phi\mapsto \langle X,\phi \rangle$ is called a 
  {\em generalised stochastic process.} This mapping is linear and continuous  
  with the usual topologies in $\cal D$ and in the space of all random variables
  $L^0(\Omega)$. Note that we can safely overload the meaning of the symbol $X$.

  The {\em mean functional\/} and the {\em correlation functional\/} of a random distribution 
  are the deterministic distribution
  $\phi\mapsto \E[\langle X,\phi\rangle]$ and the bilinear form 
  $(\phi,\psi)\mapsto \E[\langle X,\phi\rangle\langle X,\psi\rangle]$,
  respectively, provided they exist.

  A simple example of random distribution is {\em white noise\/} $\xi$, 
  characterised by the fact that
  $\langle \xi,\phi\rangle$ is centred Gaussian, with correlation functional
  $\E[\langle \xi,\phi\rangle\langle \xi,\psi\rangle]=
  \int_U \phi(s)\psi(s)\,ds$. In particular, $\langle \xi,\phi\rangle$ and
  $\langle \xi,\psi\rangle$ are independent if the supports of $\phi$ and
  $\psi$ are disjoint. 
  In this paper we will use as the base set the open half-line 
   $U=]0,+\infty[$. White noise on $U$ 
  coincides with the Wiener integral with respect to a Brownian motion $W$:
  Indeed, if $\phi$ is a test function, then 
\begin{equation}
\label{WienerInt}
  \langle \xi,\phi\rangle 
  = 
  \int_0^\infty \phi(t)\,dW(t)
\end{equation}  
  in the sense of equality in law.
  More precisely, the Wiener integral is defined as the 
  extension to $L^2(\R^+)$ of white noise (see Kuo \cite{MR0461643} for 
  a construction of the Wiener integral as extension of white noise).
  Now, integrating by parts in (\ref{WienerInt}), we can write 
\begin{equation*}
  \langle \xi,\phi\rangle 
  = 
  -\int_0^\infty W(t)\dot\phi(t)\,dt
  =
  -\langle W,\dot\phi\rangle
  \ ,
\end{equation*}  
  so that $\xi$ is the derivative of the Brownian motion $W$ as 
  random distributions. 
  A random distribution is {\em Gaussian\/} if every finite-dimensional projection
  is a Gaussian random vector. This is the case of white noise and Brownian motion.
  
  Further results on random distributions and generalised stochastic processes can
  be found for instance in the classical papers by Dawson \cite{MR41:6332} and 
  Fernique \cite{MR0221576}.
  We will also use in Section 3 the following facts about deterministic distributions, which
  apply as well to random distributions.
  
  The hyperplane $\cal H$ of $\cal D$ consisting of those functions whose integral on $U$ is equal to zero
  coincides 
  with the set of test functions which are derivatives of other test functions.
  Therefore, fixing a test function $\phi_0\in \cal D$ such that $\int_U \phi_0(t)\,dt =1$, every  
  $\phi\in \cal D$ can be uniquely decomposed as $\phi =\lambda\phi_0+\dot\psi$,
  for some $\psi\in \cal D$ and $\lambda=\int_U \phi(t)\,dt$. 
  
  If $f\in \cal D'$ is a distribution, the equation $\dot T=f$ has infinite solutions  
  (the primitives of $f$): 
  $T$ is completely determined on $\cal H$ by $\langle T,\dot\psi\rangle=-\langle f,\psi\rangle$
  whereas $\langle T,\phi_0\rangle$ can be arbitrarily chosen (for more details see
  Schwartz \cite{MR0209834}, II.4).
          

\setcounter{equation}{0}  
\section{The generalised process solution}  
\label{eu} Consider the linear \emph{stochastic differential-algebraic equation}  
(SDAE) with constant coefficients
\begin{equation}  
\label{SDAE1}  
  A\dot x + B x=f+\Lambda\xi  
  \ ,  
\end{equation}  
  where $A$ and $B$ are $n\times n$ real matrices,
  $f$ is an $n$-dimensional distribution, 
  $\Lambda$ is an $n\times m$ constant matrix, and $\xi$ is an  
  $m$-dimensional white noise: $\xi=(\xi_1,\dots,\xi_m)$, with $\xi_i$ independent
  one-dimensional white noises. 
  Recall that we will always take $U=]0,+\infty[$ as the base set for all distributions.
\begin{de}  
\label{solSDAE}
  A solution to the SDAE 
\begin{equation}\label{linSDAEcc}  
  A\dot x + Bx=f+\Lambda\xi
\end{equation}  
  is an $n$-dimensional random distribution $x$ such that, for almost all $\omega\in\Omega$,
  $x(\omega)$ is a solution to the deterministic equation
\begin{equation}\label{linSDAEcc(omega)}  
  A\dot x(\omega) + Bx(\omega)=f+\Lambda\xi(\omega)
  \ ,
\end{equation}  
  in the sense of Definition \ref{solDAE}.
\end{de}

\begin{thm}
\label{ExUnOriginalS}
  Assume $(A,B)$ is a regular matrix pencil. 
  Then, for every $u^0=(u^0_1,\dots,u^0_d)\in\R^d$, and every fixed test function
  $\phi_0$, there exists an almost surely unique random distribution $x$, 
  solution to (\ref{linSDAEcc}), 
  such that
\begin{equation*}  
  \langle Q^{-1}x,\phi_0\rangle_i=u^0_i
  \ ,
  \quad
  i=1,\dots,d
  \ .
\end{equation*}
  where $Q$ is the matrix in the reduction to KCF. Furthermore, the solution 
  is measurable with respect to the $\sigma$-field generated by $\xi$.
\end{thm}
\noindent
\emph{Proof:}
  For every $\omega\in\Omega$, we have a linear DAE with constant coefficients,
  given by (\ref{linSDAEcc(omega)}), and we know from Corollary \ref{ExUnOriginal}
  that there exists a unique solution $x(\omega)\in{\cal D}'$, satisfying 
\begin{equation*}  
  \langle Q^{-1}x(\omega),\phi_0\rangle_i=u^0_i
  \ ,
  \quad
  i=1,\dots,d
  \ .
\end{equation*}
  The fact that $x(\omega)$ is a linear and
  continuous mapping from $\cal D$ into $\R$ for every $\omega\in\Omega$ implies that $\omega\mapsto x(\omega)$
  is measurable from $(\Omega,\cal F)$ into $(\cal D',\cal B(\cal D'))$, hence a random distribution
  (see Fernique \cite{MR0221576}, section III.4). 
  We still want to prove that the mapping $\omega\mapsto x(\omega)$ is measurable
  with respect to the $\sigma$-field generated by the white noise $\xi$. 
  To this end, we will explicit the solution as much as possible with a variation    
  of constants argument.
  
%
%

  Let $P$ and $Q$ be the invertible matrices of Lemma \ref{kronecker}. Multiplying (\ref{linSDAEcc})
  from the left by $P$ and setting $y=Qx$ we obtain the SDAE in Kronecker Canonical Form
\begin{equation}  
  \label{3.2s}  
  \left(  
  \begin{matrix}  
  I_{d} & 0 \\  
  0 & N  
  \end{matrix}  
  \right)  
  \dot y  
  +  
  \left(  
  \begin{matrix}  
  J & 0 \\  
  0 & I_{q}  
  \end{matrix}  
  \right)  
  y  
  =  
  Pf+P\Lambda\xi  
  \ ,
\end{equation}  
  System (\ref{3.2s}) splits into 
  a stochastic differential system of dimension $d$
  and an  
  ``algebraic stochastic system''
  of dimension $q$.
  Denoting by $u$ and $v$ the variables in the first and the second systems
  respectively, by $b$ and $c$ the related partitioning of the 
  vector distribution $Pf$, and by $S=(\sigma_{ij})$ and $R=(\rho_{ij})$ the corresponding splitting 
  of $P\Lambda$ into matrices of dimensions $d\times m$ and $q\times m$,
  we can write the two systems as
\begin{equation}  
\label{3.3}  
  \left(  
  \begin{matrix}  
  \dot u_1 \\  
  \vdots \\  
  \dot u_d  
  \end{matrix}  
  \right)  
  +  
  J  
  \left(  
  \begin{matrix}  
  u_1 \\  
  \vdots \\  
  u_d  
  \end{matrix}  
  \right)  
  =  
  \left(  
  \begin{matrix}  
  b_1 \\  
  \vdots \\  
  b_d  
  \end{matrix}  
  \right)  
  +  
  S   
  \left(  
  \begin{matrix}  
  \xi_1 \\  
  \vdots \\  
  \xi_m  
  \end{matrix}  
  \right)  
  \ ,
\end{equation}  
\begin{equation}  
\label{3.4}  
  N  
  \left(  
  \begin{matrix}  
  \dot v_{1} \\  
  \vdots \\  
  \dot v_{q}  
  \end{matrix}  
  \right)  
  +  
  \left(  
  \begin{matrix}  
  v_{1} \\  
  \vdots \\  
  v_{q}  
  \end{matrix}  
  \right)  
  =  
  \left(  
  \begin{matrix}  
  c_{1} \\  
  \vdots \\  
  c_{q}  
  \end{matrix}  
  \right)  
  +  
  R  
  \left(  
  \begin{matrix}  
  \xi_1 \\  
  \vdots \\  
  \xi_m  
  \end{matrix}  
  \right)  
  \ .
\end{equation}  
\par  
  Fixing a test function $\phi_0$ with $\int_0^\infty \phi_0=1$ and a vector $u^0\in\R^d$,
  we have for the first one the distributional stochastic initial value problem
\begin{equation}
\label{DiffPart1}
\left.
\begin{aligned}
         \dot u +Ju &= \eta
         \\
         \langle u,\phi_0\rangle&=u^0
\end{aligned}
\right\}
\end{equation}
where $\eta:=b+S \xi$. 
  
Consider the matrix system 
\begin{equation}
\label{FundSys}
\left.
\begin{aligned}
         \dot\Phi +J\Phi &= 0
         \\
         \textstyle\int_0^\infty\Phi(t)\cdot\phi_0(t)\, dt &= I_d
\end{aligned}
\right\}
\ ,
\end{equation}
  whose distributional solution exists and is unique, and it is a ${\cal C}^{\infty}$
  matrix function $\Phi\colon\R\rightarrow\R^{d\times d}$ (see Schwartz \cite{MR0209834}, V.6). 
  Define $T:=\Phi^{-1}u$.
  From (\ref{DiffPart1}), it follows that $\dot T = \Phi^{-1}\eta$. 
  Let 
\begin{equation*}  
\Phi_{ij}(t)\phi(t)=\lambda_{ij}\phi_0(t)+\dot\psi_{ij}(t) 
\end{equation*}
be the unique decomposition of the function $\Phi_{ij}\cdot\phi\in\cal D$ into a multiple
of $\phi_0$ and an 
element of the hyperplane of derivatives $\cal H$ (see Subsection \ref{GenProc}).

Then,
\begin{align}  
  \nonumber
  \langle u_i,\phi\rangle
  &=
  \Big\langle \sum_{j=1}^d \Phi_{ij} T_j,\phi\Big\rangle
  =
  \Big\langle \sum_{j=1}^d T_j,\Phi_{ij} \phi\Big\rangle
  =
  \sum_{j=1}^d \big[ 
  \lambda_{ij} \langle T_j,\phi_0 \rangle 
  -
  \langle (\Phi^{-1}\eta)_j,\psi_{ij} \rangle \big]
  \\
  \label{gensol}
  &=
  \sum_{j=1}^d \Big[ 
   \lambda_{ij} \langle T_j,\phi_0 \rangle 
   -
   \Big\langle \sum_{k=1}^d \Phi_{jk}^{-1}\eta_k,\psi_{ij} \Big\rangle \Big]
  =
  \sum_{j=1}^d  
   \lambda_{ij} \langle T_j,\phi_0 \rangle 
   -
    \sum_{k=1}^d \Big\langle\eta_k,\sum_{j=1}^d\psi_{ij}\Phi_{jk}^{-1} \Big\rangle 
  \ .  
\end{align}  
  The terms $\langle T_j,\phi_0 \rangle$ should be defined in order to fulfil 
  the initial condition. 
  Using the decomposition 
\begin{equation*}  
\Phi_{ij}(t)\phi_0(t)=\delta_{ij}\phi_0(t)+\dot\psi^0_{ij}(t) 
\end{equation*}
  and applying
  formula (\ref{gensol}) to $\phi=\phi_0$, 
  it is easily found that we must define  
\begin{equation*}  
  \langle T_j,\phi_0\rangle
  =
  u_j^0 + \sum_{k=1}^d \Big\langle \eta_k,\sum_{\ell=1}^d \psi_{j\ell}^0\Phi_{\ell k}^{-1}\Big\rangle
  \ .
\end{equation*}
  Therefore,   
\begin{align*}  
  \langle u_i,\phi\rangle
  &=
  \sum_{j=1}^d  
   \lambda_{ij} 
    \Big[
     u^0_j + \sum_{k=1}^d \Big\langle \eta_k,\sum_{\ell=1}^d\psi_{j\ell}^0\Phi_{\ell k}^{-1} \Big\rangle 
    \Big]
   -
    \sum_{k=1}^d \Big\langle\eta_k,\sum_{j=1}^d\psi_{ij}\Phi_{jk}^{-1} \Big\rangle 
  \\  
  &=
  \sum_{j=1}^d  
    \lambda_{ij} 
    u^0_j 
   +
    \sum_{k=1}^d \Big\langle \eta_k,\sum_{j=1}^d\lambda_{ij}\sum_{\ell=1}^d\psi_{j\ell}^0\Phi_{\ell k}^{-1} 
    -\sum_{j=1}^d\psi_{ij}\Phi_{jk}^{-1} \Big\rangle
  \\  
  &=
  \sum_{j=1}^d  
    \lambda_{ij} 
    u^0_j 
   +
    \sum_{k=1}^d \Big\langle \eta_k,\sum_{\ell=1}^d
    \Big(    
    \sum_{j=1}^d\lambda_{ij}\psi_{j\ell}^0-\psi_{i\ell}\Big)\Phi_{\ell k}^{-1} 
    \Big\rangle
  \ .  
\end{align*}  
  Taking into account that 
\begin{align*}  
  \psi_{j\ell}^0 (t)
  &=
  \int_0^t
  \big(
  \Phi_{j\ell}(s) \phi_0(s)
  -
  \delta_{j\ell} \phi_0(s)\big)
  \, ds
  \ ,
  \\
  \psi_{i\ell} (t)
  &=
  \int_0^t
  \big(
  \Phi_{i\ell}(s) \phi(s)
  -
  \lambda_{i\ell} \phi_0(s)\big)
  \, ds
  \ ,
\end{align*}
 we obtain finally
\begin{align}\label{crm1}
  \nonumber
  \langle u_i,\phi\rangle
  &=
  \sum_{j=1}^d  
    \lambda_{ij} 
    u^0_j 
   +
    \sum_{k=1}^d \Big\langle \eta_k,\sum_{\ell=1}^d
    \Big(    
    \int_0^t
    \Big(
    \sum_{j=1}^d\lambda_{ij}\Phi_{j\ell}(s)\phi_0(s)-\Phi_{i\ell}(s)\phi(s)\Big)\,ds\Big)\Phi_{\ell k}^{-1}(t) 
    \Big\rangle
  \\   
  \nonumber
  &=
  \sum_{j=1}^d  
    \lambda_{ij} 
    u^0_j 
   +
    \sum_{k=1}^d \Big\langle \eta_k,\sum_{\ell=1}^d
    \Big(    
    \int_0^t
    (
    \lambda\Phi\phi_0-\Phi\phi)(s)\,ds\Big)_{i\ell}
    \Phi_{\ell k}^{-1}(t) 
    \Big\rangle
  \\  
  &=
  \sum_{j=1}^d  
    \lambda_{ij} 
    u^0_j 
   +
    \sum_{k=1}^d \Big\langle \eta_k,
    \Big(    
    \int_0^t
    (
    \lambda\Phi\phi_0-\Phi\phi)(s)\,ds
    \cdot
    \Phi^{-1}(t)\Big)_{ik} 
    \Big\rangle
  \ .  
\end{align}  

  
  On the other hand, the algebraic part (\ref{3.4}) consists of a number of decoupled blocks, 
  which are easily solved by backward substitution. 
  Any given block can be solved independently of the others and  
  a recursive calculation gives, e.g. for a first block of dimension $q_1$, 
  the following generalised process  
  solution
\begin{equation}\label{sae.dist}
  \langle v_j,\phi\rangle
  =
  \sum_{k=j}^{q_1}  
  \Big\langle 
  c_k+\sum_{\ell=1}^{m} \rho_{k\ell}\xi_\ell,\phi^{(k-j)}
  \Big\rangle
  \ ,
  \quad j=1,\dots q_1 \ ,
  \quad \phi\in {\cal D} \ .
\end{equation}  
  By (\ref{crm1}) and (\ref{sae.dist}), we have $(u,v)=G(\xi)$, for some deterministic
  function $G\colon ({\cal D'})^m \rightarrow ({\cal D'})^n$.
  Given generalized sequence $\{\eta_\alpha\}_\alpha \subset ({\cal D'})^m $ converging
  to $\eta$ in the product of weak-$\star$ topologies, it is immediate to see that
  $G(\eta_\alpha)$ converges to $G(\eta)$, again in the product of weak-$\star$ topologies.
  This implies that the mapping $G$ is continuous and therefore
  measurable with respect to the Borel $\sigma$-fields. Thus, the solution process $x$ 
  is measurable with respect to the $\sigma$-field generated by $\xi$.
\qed  
  
\begin{re}  
\label{re3.3}
  In the case when $b=0$, so that the right hand side in (\ref{DiffPart1}) is simply $S\xi$, 
  it is well known that the solution of the differential system is a classical stochastic process 
  which can be expressed as a functional of a standard $m$-dimensional Wiener process. Indeed, 
  we have, in the sense of equality in law, from (\ref{crm1}),
\begin{equation*}
  \langle u_i,\phi\rangle
  =
  \sum_{j=1}^d \lambda_{ij}u_j^0
  +
  \sum_{k=1}^d   \sum_{\ell=1}^m 
  \int_0^{\infty} 
  \Big(
   \int_0^t 
   \big( 
    \lambda\Phi\phi_0-\Phi\phi
    \big)
   (s)\, ds 
   \cdot 
   \Phi^{-1}(t)
   \Big)_{ik}
  \sigma_{k\ell} 
  \, dW_\ell(t)
  \ .
\end{equation*}
  Fix an initial time $t_0\in]0,\infty[$. Take a sequence $\{\phi_0^n\}_n\subset \cal D$ converging
  in $\cal D'$ to the Dirac delta 
  $\delta_{t_0}$, and with $\supp\, \phi_0^n\subset[t_0-\frac{1}{n},t_0+\frac{1}{n}]$, 
  and let $\{\Phi^n\}_n$ be the corresponding sequence of solutions to the matrix system
  (\ref{FundSys}).
  Then,     
  $\lim_{n\to\infty} \int_0^t \Phi^n\phi_0^n= I_d\cdot\mathbf 1_{[t_0,\infty[}(t)$ a.e. 
  and we get
\begin{align*}
  \langle u_i,\phi\rangle
  =
  \sum_{j=1}^d \lambda_{ij}u_j^0
  &+
  \sum_{k=1}^d   \sum_{\ell=1}^m 
  \int_{t_0}^{\infty} 
  \big(
   \lambda\Phi^{-1}\big)_{ik}(s) \sigma_{k\ell}\,dW_\ell(s) 
  \\ 
  &-
  \sum_{k=1}^d   \sum_{\ell=1}^m 
  \int_{0}^{\infty} 
  \Big(
  \int_{0}^{s} 
  \big(
   \Phi\cdot\phi\big)(u)\,du\cdot\Phi^{-1}(s)\Big)_{ik} \sigma_{k\ell}\,dW_\ell(s) 
   \ ,  
\end{align*}
  where $\Phi$ and $\lambda$ are now the limit of their corresponding sequences.  
  
  Now collapsing in the same way $\phi$ to $\delta_t$, which implies the convergence of $\lambda$
  to $\Phi$ and that of 
  $\int_0^s \Phi\phi$ to $I_d\cdot\mathbf 1_{[t,\infty[}(s)$ a.e.,
  we arrive at 
\begin{equation*}
  u_i(t)
  =
  \sum_{j=1}^d \Phi_{ij}u_j^0
  +
  \sum_{k=1}^d   \sum_{\ell=1}^m 
  \int_{t_0}^{t} 
  \big(
   \Phi(t)\Phi^{-1}(s)\big)_{ik} \sigma_{k\ell}\,dW_\ell(s) 
   \ .
\end{equation*}
  Finally, using that the solution to (\ref{FundSys}) with $\delta_{t_0}$ in place of $\phi_0$
  is known to be $\Phi(t)=e^{-J(t-t_0)}$, we obtain
\begin{equation*}
  u(t)
  =
  e^{-J(t-t_0)}  
  \Big[
  u_0
  +
  \int_{t_0}^{t} 
  e^{-J(s-t_0)}S\,dW(s) 
  \Big]
   \ .
\end{equation*}
In a similar way we can express the first block of the algebraic part, if $c=0$, as
\begin{equation}
  \langle v_j,\phi\rangle
  =
  \sum_{k=j}^{q_1}  
  \sum_{\ell=1}^m  
  \rho_{k\ell}
  \int_0^\infty  
  \phi^{(k-j)}(t)\,dW_\ell(t)
  \ ,
  \quad j=1,\dots q_1
  \ ,
\end{equation}  
  and analogously for any other block.
\qed
\end{re}

%
%
%

\setcounter{equation}{0}  
\section{The law of the solution}  
\label{ac}  
  In the previous section we have seen that the solution to a linear SDAE
  with regular pencil and additive white noise
  can be explicitly given as a functional of the input noise.
  From the modelling viewpoint, the law of the solution is the important
  output of the model. Using the explicit form of the solution, one can try to investigate
  the features of the law in which one might be interested. 

  To illustrate this point, we shall study the absolute continuity properties of the 
  joint law of the vector solution  
  evaluated at a fixed arbitrary test function $\phi$.
  We will assume throughout this section that the base probability space is the canonical
  space of white noise: $\Omega={\cal D}'$, ${\cal F}={\cal B}({\cal D}')$, and 
  $P$ is the law of white noise. This will be used in Theorem
  \ref{Thm43}, to ensure the existence of conditional probabilities (see 
  Dawson \cite{MR41:6332}, Theorem 2.12).
  
  Let us start by considering separately the solutions to the decoupled equations
  (\ref{3.3}) and (\ref{3.4}). From the explicit calculation in the
  previous section (equation (\ref{crm1}) for the differential part and 
  equation (\ref{sae.dist}) for the first algebraic block), we get that
  for any given test function $\phi$ the random vectors $\langle u,\phi \rangle$  
  and $\langle v,\phi \rangle$ 
  have a Gaussian distribution with expectations
\begin{align*}
  E[\langle u_i,\phi \rangle] 
  &=  
  \sum_{j=1}^d  
    \lambda_{ij} 
    u^0_j 
   +
    \sum_{k=1}^d \Big\langle b_k,
    \Big(    
    \int_0^t
    (
    \lambda\Phi\phi_0-\Phi\phi)(s)\,ds
    \cdot
    \Phi^{-1}(t)\Big)_{ik} 
    \Big\rangle
  \ ,\quad
  i=1,\ldots,d
  \ ,
  \\
  E[\langle v_i,\phi \rangle] 
  &=
  \sum_{k=i}^{q_1}  
  \langle 
  c_k,\phi^{(k-i)}
  \rangle
  \ ,\quad
  i=1,\ldots,q_1
  \ ,
\end{align*}
  and covariances
\begin{align}
\label{covD} 
\Cov\big(\langle u,\phi \rangle\big)_{ij} 
  &=
  \sum_{\ell=1}^m 
  \int_0^\infty  
  \Big[ 
    \sum_{k=1}^d     
    \int_0^t
    (\lambda\Phi\phi_0-\Phi\phi)(s)\,ds
    \cdot
    \Phi^{-1}(t))_{ik} \sigma_{kl}\Big]^2
  \\ 
  \notag
  &\phantom{\ \sum_{k=1}^d}     
  \times
  \Big[ 
    \sum_{k=1}^d     
    \int_0^t
    (\lambda\Phi\phi_0-\Phi\phi)(s)\,ds
    \cdot
    \Phi^{-1}(t))_{jk} \sigma_{kl}\Big]^2
  \,ds
  \ ,\quad
  i,j=1,\ldots,d
  \ ,
\end{align}
and
\begin{align}
\label{covalgblock}
   \Cov \big(\langle v,\phi\rangle\big)
  =&
   \left(
   \begin{matrix}   
     \rho_1 & \rho_2 & \cdots & \rho_{q_1-1} & \rho_{q_1}
     \\
     \rho_2 & \rho_3 & \cdots & \rho_{q_1} & 0
     \\
     \vdots & & & & \vdots 
     \\
     \rho_{q_1} & 0 & \cdots & 0 & 0
   \end{matrix}
   \right)
  \nonumber
   \\
   &  \times
   \Cov
   \Big(
   \langle \xi,\phi\rangle 
   ,
   \langle \xi,\dot\phi\rangle 
   ,
   \dots
   ,
   \langle \xi,\phi^{(q_1-1)}\rangle 
   \Big)
   \left(
   \begin{matrix}   
     \rho_1 & \rho_2 & \cdots & \rho_{q_1-1} & \rho_{q_1}
     \\
     \rho_2 & \rho_3 & \cdots & \rho_{q_1} & 0
     \\
     \vdots & & & & \vdots 
     \\
     \rho_{q_1} & 0 & \cdots & 0 & 0
   \end{matrix}
   \right)^{\top}
   \ ,
\end{align}
  where $\rho_i$ denotes the $i$-th row of the matrix $R$ and
  $ \Cov\big(\langle \xi,\phi\rangle, \dots,
   \langle \xi,\phi^{(q_1-1)}\rangle \big)$ is a square matrix
  of dimension $m q_1$. 
  We refer the reader to \cite{MR1956867} for a comprehensive 
  study of multidimensional Gaussian laws.

  For the differential variables $u$ alone, we are faced with a usual 
  linear stochastic differential equation (see Remark \ref{re3.3}), and there are 
  well-known results on sufficient conditions for its absolute continuity, involving   
  the matrices $S$ and $J$ (see e.g. Nualart \cite{MR2200233}). 
  
  For the algebraic variables $v$, their absolute continuity depends in part on the
  invertibility of the covariance matrix of the white noise and its derivatives
  that appear in (\ref{covalgblock}). We will use the following auxiliary result
  concerning the joint distribution of 
  a one-dimensional white noise and its first $k$ derivatives.
  This is a vector distribution with a centred Gaussian law and a covariance 
  that can be expressed in full generality as (cf. Subsection 2.2)
\begin{equation}\label{cov}  
  \Cov\big(\langle \xi,\phi\rangle, \dots, \langle \xi^{(k)},\phi\rangle\big)_{ij}
  =
  \realpart \Big[(-1)^{\frac{|i-j|}{2}}\Big]\| \phi^{((i+j)/2)} \|^2
  \ ,
\end{equation}
  where $\realpart$ means the real part. We can prove the absolute continuity
  of this vector for $k\le 3$.
\begin{lm}  
\label{AbsContWN1}
  For all $\phi\in{\cal D}-\{0\}$,
  and a one-dimensional white noise $\xi$, the vector $\langle(\xi,\dot\xi,\ddot\xi,\dddot\xi),\phi\rangle$
  is absolutely continuous.
\end{lm}
{\it Proof:}
  The covariance matrix of the vector $\langle(\xi,\dot\xi,\ddot\xi,\dddot\xi),\phi\rangle$
  is
$$
\left(
\begin{matrix}
  \|\phi\|^2 & 0 & -\|\dot\phi\|^2 & 0
  \\
  0 & \|\dot\phi\|^2 & 0  & -\|\ddot\phi\|^2
  \\
  -\|\dot\phi\|^2 & 0  & \|\ddot\phi\|^2 & 0
  \\
  0 & -\|\ddot\phi\|^2 & 0 & \|\dddot\phi\|^2
\end{matrix}
\right)
$$
  whose determinant is equal to 
$$
  \det
  \left(
  \begin{matrix}
    \|\phi\|^2 & -\|\dot\phi\|^2
    \\
    -\|\dot\phi\|^2 & \|\ddot\phi\|^2
  \end{matrix}
  \right)
  \cdot
  \det
  \left(
  \begin{matrix}
  \|\dot\phi\|^2 & -\|\ddot\phi\|^2
  \\
  -\|\ddot\phi\|^2 & \|\dddot\phi\|^2
\end{matrix}
\right)\ .
$$
  Both factors are strictly positive, in view of the chain of strict inequalities 
\begin{equation}
\label{cs1}
  \frac{\|\phi\|}{\|\dot\phi\|}
  >
  \frac{\|\dot\phi\|}{\|\ddot\phi\|}
  >
  \frac{\|\ddot\phi\|}{\raise-1.7pt\hbox{$\|\dddot\phi\|$}}
  >
  \cdots
  \ ,
  \qquad
  \phi\in{\cal D},\ \phi\not\equiv 0
\end{equation}
  These follow
  from integration by parts and Cauchy-Schwarz inequality, e.g.
$$
  \|\dot\phi\|^2
  =
  \int_0^{\infty} \dot\phi\cdot\dot\phi
  =
  -\int_0^{\infty} \phi\cdot\ddot\phi
  \le
  \|\phi\|\cdot\|\ddot\phi\|
  \ , 
  $$
  and the inequality is strict unless $\ddot\phi=K\phi$ for some $K$, which 
  implies $\phi\equiv 0$.
\qed

\noindent
  The proof above does not longer work for higher order derivatives and we do not know
  if the result is true or false.

  Consider, as in the previous section, only the first algebraic block, and 
  assume momentarily that its dimension is $q_1=2$. 
  From (\ref{covalgblock}), the covariance matrix of the random vector 
  $\langle (v_1,v_2),\phi\rangle$ is 
\begin{equation*} 
   \left(
   \begin{matrix}   
     \|\phi\|^2\|\rho_1\|^2+ \|\dot\phi\|^2\|\rho_2\|^2 
    &
     \|\phi\|^2\langle\rho_1,\rho_2\rangle 
    \vspace{2mm}\\ 
     \|\phi\|^2\langle\rho_1,\rho_2\rangle 
    & 
     \|\phi\|^2\|\rho_2\|^2 
   \end{matrix} 
   \right)
   \ ,
\end{equation*} 
with determinant  
\begin{equation*}
  \|\phi\|^4
  \big(\|\rho_1\|^2 \|\rho_2\|^2 - \langle \rho_1,\rho_2\rangle^2\big)
  +
  \|\phi\|^2\|\dot\phi\|^2\|\rho_2\|^4
  \ .
\end{equation*}    
  Hence, assuming $\phi\not\equiv0$, we see that  the joint law of $\langle v_1,\phi\rangle$ and $\langle v_2,\phi\rangle$
  is absolutely continuous with respect to Lebesgue measure in $\R^2$ if $\rho_2$ is not
  the zero vector. 
  When $\|\rho_2\|=0$ but $\|\rho_1\|\neq 0$, then $\langle v_2,\phi\rangle$ is degenerate
  and $\langle v_1,\phi\rangle$ is absolutely continuous, whereas 
  $\|\rho_2\|=\|\rho_1\|=0$ makes the joint law degenerate to a point.
  
  This sort of elementary analysis, with validity for any test function $\phi$, can be carried out 
  for algebraic blocks of any nilpotency index, as it is proved in the next proposition. 
  Let us denote by ${\cal E}(k)$ the subset of test functions $\phi$ such that the covariance 
  $\Cov\big(\langle \xi,\phi\rangle,\dots,\langle \xi^{(k-1)},\phi\rangle\big)$ is nonsingular.
  With an $m$-dimensional white noise, the covariance is a matrix with
  $(k+1)^2$ square 
  $m\times m$ blocks, where the block $(i,j)$ is
  $\realpart \Big[(-1)^{\frac{|i-j|}{2}}\Big]\| \phi^{((i+j)/2)} \|^2$ times the identity $I_m$.
\begin{pr}
\label{AlgPart}
  Let $(v_1,\dots,v_{q_1})$ be the generalised process solution to the first block of the algebraic system
  (\ref{3.4new}) and $r$ the greatest row index such that $\|\rho_r\|\neq 0$, and fix
  $\phi\in {\cal E}(q_1)$.  

\noindent
  Then $\langle (v_1,\dots,v_r),\phi\rangle$ is a Gaussian absolutely continuous
  random vector and $\langle (v_{r+1},\dots,v_{q_1}),\phi\rangle$ degenerates
  to a point.
\end{pr}
\noindent
\emph{Proof:} We can assume that
  $c=0$, since the terms $\sum_{k=j}^{q_1} \langle c_k,\phi^{(k-j)}\rangle$
  in (\ref{3.4new}) only contribute as additive constants. Then we can write
\[
  \left(  
  \begin{matrix}  
  \langle v_1,\phi\rangle \\  
  \vdots \\  
  \langle v_{q_1},\phi\rangle 
  \end{matrix}  
  \right)  
  =
   \left(
   \begin{matrix}   
     \rho_1 & \rho_2 & \cdots & \rho_{q_1-1} & \rho_{q_1}
     \\
     \rho_2 & \rho_3 & \cdots & \rho_{q_1} & 0
     \\
     \vdots & & & & \vdots   
     \\
     \rho_{q_1} & 0 & \cdots & 0 & 0
   \end{matrix}
   \right)
    \left(
   \begin{matrix}   
     \langle \xi,\phi\rangle     \\
     \langle \xi,\dot \phi\rangle
     \\
     \vdots
     \\
     \langle \xi,\phi^{(q_1-1)}\rangle
   \end{matrix}
   \right)
   \ .
\]
  If $r$ is the greatest row index with $\|\rho_r\|\neq 0$, it is clear that the
  $q_1 \times m q_1$ matrix 
\[
   \left(
   \begin{matrix}   
     \rho_1 & \rho_2 & \cdots & \rho_{q_1-1} & \rho_{q_1}
     \\
     \rho_2 & \rho_3 & \cdots & \rho_{q_1} & 0
     \\
     \vdots & & 
     \\
     \rho_{q_1} & 0 & \cdots & 0 & 0
   \end{matrix}
   \right)
   \ ,
\]
  has rank $r$. The linear transformation  
  given by this matrix is onto $\R^{r}\times \{0\}^{q_1-r}$. 
  From this fact and the absolute continuity of the vector
  $(\langle \xi,\phi\rangle,\dots,\langle \xi^{(q_1-1)},\phi\rangle)$, it
  is immediate that the vector $(\langle v_1,\phi\rangle,\dots,\langle v_r,\phi\rangle)$
  is absolutely continuous, while $(\langle v_{r+1},\phi\rangle,\dots,\langle v_{q_1},\phi\rangle)$
  degenerates to a point. \qed

  Let us now consider the solution $x$ to the whole SDAE (\ref{linSDAEcc}).
  We will state a sufficient condition for the absolute continuity of
  $\langle x,\phi\rangle$, $\phi\in \cal D$.
  The following
  standard result in linear algebra will be used (see e.g. Horn and
  Johnson \cite{Horn}, page 21).
\begin{lm}
\label{schur}
  Let the real matrix M read blockwise
\[
M =
   \left(
   \begin{matrix}   
     A & B
     \\
     C & D
   \end{matrix}
   \right)
\]
  where $A \in \R^{d \times d}$, $B \in \R^{d \times q}$,  
  $C \in \R^{q \times d}$, $D \in \R^{q \times q}$ and $D$ is
  invertible. Then the $d\times d$ matrix
\[
A - B D^{-1}C
\]
  is called the \emph{Schur complement} of $D$ in $M$ and it holds that
\[
\det M=\det D \cdot \det (A - B D^{-1}C) 
\ .
\]
\end{lm}
\qed

\noindent
  A natural application of this lemma is in solving a system of
  linear equations:
\[
\left.
\begin{aligned}
         Ax+By =u
         \\
         Cx+Dy=v
\end{aligned}
\right\}
\]
  where $x,u\in \R^d$, $y,v\in\R^q$.
  We have
\begin{equation}
\label{schurSolve}
u = (A - B D^{-1}C)x +BD^{-1}v 
\end{equation}
  and, if $M$ is in addition invertible, the solution to the linear equation
  is given by
\begin{align*}
         x & =(A - B D^{-1}C)^{-1} (u - BD^{-1}v)
         \\
         y & =D^{-1} (v-C(A - B D^{-1}C)^{-1}(u-BD^{-1})) \ .
\end{align*}

  We now state and prove the main result of this section.
\begin{thm}
\label{Thm43}
  Assume $(A,B)$ is a regular matrix pencil and that the matrix
  $\Lambda$ of equation (\ref{SDAE1}) has full rank, and
  call $r$ the nilpotency index of the SDAE (\ref{SDAE1}).
  Then the law of the unique solution to the SDAE (\ref{SDAE1})
  at any test function $\phi\in {\cal E}(r)$ is absolutely continuous
  with respect to Lebesgue measure on $\R^n$.
\end{thm}

\begin{co}
\label{cor44}
  Under the assumptions of Theorem \ref{Thm43}, if the nilpotency index is $r \leq 4$,
  then the law is absolutely continuous at every test function 
  $\phi\in {\cal D}-\{0\}$. 
\end{co}
{\it Proof of Theorem \ref{Thm43}:}
  It will be enough to prove that the random vector
  $\langle (u,v),\phi\rangle $, solution to (\ref{3.3}) and (\ref{3.4}),  
  admits an absolutely
  continuous law at any test function $\phi\in \cal E$,
  since the solution to the original system is then obtained through the non-singular
  transformation $Q$.

  We shall proceed in two steps: first we shall prove that $\langle v,\phi\rangle $ 
  admits an absolutely continuous law, and then that the conditional law of 
  $\langle u,\phi\rangle $, given 
  $\langle v,\phi\rangle $, is also absolutely continuous, almost surely with respect to 
  the law of $\langle v,\phi\rangle $.

  \emph{Step 1:} We can assume $c=0$ in (\ref{3.4}). By Proposition \ref{AlgPart}, the solution to any
  algebraic block is separately absolutely continuous.
  Assume now that
  there are exactly two blocks of dimensions $q_1$ and $q_2$, with $q_2\le q_1$; 
  the case with an arbitrary number
  of blocks does not pose additional difficulties. 
%

 As in Proposition \ref{AlgPart},
  we have 
\begin{align}
  \label{bigmatrix}
  \left(  
  \langle v_1,\phi\rangle, \ldots, \langle v_{q_1},\phi\rangle, 
  \right.
  & 
  \left. \langle v_{q_1+1},\phi\rangle,
  \ldots, \langle v_{q_1+q_2},\phi\rangle
  \right)^{\top}
  =
  \nonumber
  \\
  & =   
  \left(
   \begin{matrix}   
     \rho_1 & \rho_2 & \cdots & \cdots & \rho_{q_1-1} & \rho_{q_1}
     \\
     \rho_2 & \rho_3 & \cdots & \cdots & \rho_{q_1} & 0
     \\
     \vdots & & & & & \vdots   
     \\
     \rho_{q_1} & 0 & \cdots & 0 & 0 & 0
     \\
     \rho_{q_1+1} & \rho_{q_1+2} & \cdots & \cdots & \rho_{q_1+q_2} & 0
     \\
     \rho_{q_1+2} & \rho_{q_1+3} & \cdots & \rho_{q_1+q_2} & 0 & 0
     \\
     \vdots & & & & & \vdots   
     \\
     \rho_{q_1+q_2} & 0 & \cdots & 0 & 0 & 0
        \end{matrix}
   \right)
    \left(
   \begin{matrix}   
     \langle \xi,\phi\rangle     \\
     \langle \xi,\dot \phi\rangle
     \\
     \vdots
     \\
     \langle \xi,\phi^{(q_1-1)}\rangle
   \end{matrix}
   \right)
\end{align}
  Since the $(q_1+q_2)\times m$ matrix $R=(\rho_1,\dots,\rho_{q_1+q_2})^{\top}$ has, 
  by the hypothesis on $\Lambda$, 
  full rank equal to $q_1+q_2$,
  the transformation 
  defined by (\ref{bigmatrix}) is onto $\R^q$.
  From the absolute continuity of the vector 
  $(\langle \xi,\phi\rangle, \dots, \langle \xi^{(q_1-1)},\phi\rangle)$,
  we deduce that of $\langle v,\phi\rangle$.

  \emph{Step 2:}  
  Since the $n\times m$ matrix $\Lambda$ has full rank 
  we can assume, 
  reordering columns if necessary,
  that the submatrix made of
  its first $n$ rows and columns is invertible, and that
\[
   \Lambda=
   \left(
   \begin{matrix}   
     S
     \\
     R
   \end{matrix}
   \right) 
   =
   \left(
   \begin{matrix}   
     A & B & E
     \\
     C & D & F
   \end{matrix}
   \right) \ ,
\]
  where $A \in \R^{d\times d}$, $B \in \R^{d\times q}$, $C \in \R^{q \times d}$,   
   $D \in \R^{q\times q}$, $E \in \R^{d\times (m-n)}$, $F \in \R^{q\times (m-n)}$, 
  with invertible $D$. 
  Let us define 
\[  
  w:=(v_1, \ldots, v_q, v_1 + \dot v_2, \ldots, 
  v_{q_1-1} + \dot v_{q_1}, v_{q_1+1} + \dot v_{q_1+2}, \dots, v_{q_1+q_2-1} + \dot v_{q_1+q_2},
  \xi_{n+1}, \ldots, \xi_m)
  \ .
\]  
  We can write then
\begin{equation}
  \label{verybigmatrix}
  \langle (w,\xi_1, \ldots, \xi_d),\phi\rangle^{\top}
  =   
  \left(
   \begin{matrix}   
   G_1
   \\
   G_2
   \\
   G_3
   \\
   G_4
   \\
   G_5
   \end{matrix}
   \right)
    \left(
   \begin{matrix}   
     \langle \xi,\phi\rangle     \\
     \langle \xi,\dot \phi\rangle
     \\
     \vdots
     \\
     \langle \xi,\phi^{(q_1-1)})\rangle
   \end{matrix}
   \right) \ ,
\end{equation}
  where $G_1$ is the matrix in (\ref{bigmatrix}), 
\begin{align*}
  &G_2 :=
  \left(
   \begin{matrix}   
     \rho_{1} & 0 & \cdots & 0
     \\
     \vdots & \vdots & &  \vdots   
     \\
     \rho_{q_1-1} & 0 & \cdots & 0
   \end{matrix}
  \right)
  \ , \quad
  G_3 :=
  \left(
   \begin{matrix}   
     \rho_{q_1+1} & 0 & \cdots & 0
     \\
     \vdots & \vdots & &  \vdots   
     \\
     \rho_{q_1+q_2-1} & 0 & \cdots & 0
   \end{matrix}
  \right)
  \ ,
\\[5mm]
   & G_4 :=
 \left(
 \begin{matrix}   
     e_{d+q+1} & 0 & \cdots & 0
     \\
     \vdots & \vdots & &  \vdots   
     \\
     e_{m} & 0 & \cdots & 0
   \end{matrix}
  \right)
  \ , \quad
  G_5 :=
  \left(
   \begin{matrix}   
     e_{1} & 0 & \cdots & 0
     \\
     \vdots & \vdots & &  \vdots   
     \\
     e_{d} & 0 & \cdots & 0
   \end{matrix}
  \right)
  \ ,
\end{align*}
  and $(e_i)_j=\delta_{ij}$.

  By the invertibility of $D$ and the fact that the rows $\rho_{q_1}$ and $\rho_{q_1+q_2}$
  have at least one element different from zero,
  it is easy to see that the matrix in (\ref{verybigmatrix})
  has itself full rank and, therefore, that the vector 
  $\langle (w,\xi_1, \ldots, \xi_d),\phi\rangle$ has an absolutely continuous Gaussian law.

  Using (\ref{schurSolve}), we obtain 
\begin{equation}
\label{from4.5}
  \left(  
  \begin{matrix}  
  \dot u_1 \\  
  \vdots \\  
  \dot u_d  
  \end{matrix}  
  \right)  
  +  
  J  
  \left(  
  \begin{matrix}  
  u_1 \\  
  \vdots \\  
  u_d  
  \end{matrix}  
  \right)  
  =  
  (A - B D^{-1}C)
  \left(  
  \begin{matrix}  
  \xi_1 \\  
  \vdots \\  
  \xi_{d}  
  \end{matrix}  
  \right)  
  +
  BD^{-1} 
  \left(  
  \begin{matrix}  
  v_1 + \dot v_2 \\
  v_2 + \dot v_3 \\
  \vdots \\  
  v_{q_1} \\  
  v_{q_1+1} + \dot v_{q_1+2} \\
  v_{q_1+2} + \dot v_{q_1+3} \\
  \vdots \\  
  v_{q_1+q_2}  
  \end{matrix}  
  \right)
  + 
  (E- BD^{-1}F) 
  \left(  
  \begin{matrix}  
  \xi_{n+1} \\  
  \vdots \\  
  \xi_{m}  
  \end{matrix}  
  \right)  
  \ .
\end{equation}
  
  It is obvious that both Definition \ref{solSDAE} and Theorem \ref{ExUnOriginalS} continue to hold true 
  with any generalised process $\theta$ in place of the white
  noise $\xi$ in the right-hand side.
  From Theorem \ref{ExUnOriginalS} we have in particular that the solution $u$ 
  to the differential system (\ref{from4.5}) is a measurable
  function $G\colon({\cal D}')^m\rightarrow ({\cal D}')^d$ of its right-hand side $\theta$, that means, $u=G(\theta)$.
  Let 
\begin{equation*}  
  p:{\cal B}({\cal D}')\times{\cal D}'\longrightarrow[0,1]
  \quad
  \text{and}
  \quad
  q:{\cal B}({\cal D}')\times{\cal D}'\longrightarrow[0,1]
\end{equation*}
  be conditional laws of $u$ given $w$, and of $\theta$ given $w$, respectively. 
  That means that
\begin{equation*}  
  P\big(\{u\in B\}\cap \{w\in C\}\big)
  =
  \int_{C} p(B,w)\,\mu(dw)
\end{equation*}
and  
\begin{equation*}  
  P\big(\{\theta\in B\}\cap \{w\in C\}\big)
  =
  \int_{C} q(B,w)\,\mu(dw)
  \ ,
\end{equation*}
  for any $B,C\in{\cal B}({\cal D}')$, where $\mu$ is the law of $w$.

  For every $w\in{\cal D}'$, let $Z_{w}$ be a random distribution
  $Z_{w}\colon\Omega\rightarrow {\cal D}'$ with law
  $P\{Z_{w}\in B\}=q(B,w)$.
  
  Then 
\begin{align*}  
  &
  \int_{C} P\{G(Z_{w})\in B\}\,\mu(dw)
  =
  \int_{C} P\{Z_{w}\in G^{-1}(B)\}\,\mu(dw)
  \\
  &=
  \int_{C} q(G^{-1}(B),w)\,\mu(dw)
  =
  P\big(\{\theta\in G^{-1}(B)\}\cap\{w\in C\}\big)
  \\
  &=
  P\big(\{G(\theta)\in B\}\cap\{w\in C\}\big)
  =
  P\big(\{u\in B\}\cap\{w\in C\}\big)
  \\
  &=
  \int_{C} p(B,w)\,\mu(dw)
\end{align*}
  
  Therefore $P\{G(Z_{w})\in B\}=p(B,w)$ almost surely with respect
  to the law of $w$, for all $B\in\cal B(\cal D')$. 
  We have proved that if the right-hand side of the differential system has the
  law of $\theta$ conditioned to $w$, then its solution has the law of $u$
  conditioned to $w$. It remains to show that this conditional law 
  is absolutely continuous, almost surely with respect to the law of $w$.
  
  Now, for each $w$, we can take $Z_{w}$ as 
\begin{equation*}  
  Z_{w} =  
    (A - B D^{-1}C)
  \eta_w
  + a_w
\end{equation*}  
  where $a_w$ is a constant $d$-dimensional distribution, and $\langle\eta,\phi\rangle$ is,   
  for each $\phi\in\cal D$, a Gaussian $d$-dimensional vector. This random
  vector is absolutely continuous: Indeed, its law is that of the $d$ first components of the 
  $m$-dimensional white noise $(\xi_1,\dots,\xi_m)$ conditioned to lie in an $(m-d)$-dimensional
  linear submanifold.  
  Let $L_{w,\phi}$ be its covariance matrix. 
  Then $\langle\eta,\phi\rangle=L_{\omega,\phi}^{1/2}\langle\zeta,\phi\rangle$, for some 
  $d$-dimensional white noise $\zeta=(\zeta_1,\dots,\zeta_d)$. 
  
  Consider now the (ordinary) stochastic differential equation
\begin{equation}
\label{newcompactform}
  \left(  
  \begin{matrix}  
  \langle\dot{\underline u}_1,\phi\rangle \\
  \vdots \\  
  \langle\dot{\underline u}_d,\phi\rangle
  \end{matrix}  
  \right) 
  + 
  J
  \left(  
  \begin{matrix}  
  \langle{\underline u}_1,\phi\rangle \\
  \vdots \\  
  \langle{\underline u}_d,\phi\rangle
  \end{matrix}  
  \right) 
  =
  a_w
  +
  (A - B D^{-1}C)
  L_{w,\phi}^{1/2}
  \left(  
  \begin{matrix}  
  \langle\zeta_1,\phi\rangle \\
  \vdots \\    
  \langle\zeta_d,\phi\rangle
  \end{matrix}  
  \right)    
  \ .
\end{equation}
  
  By hypothesis, the Schur complement $A - B D^{-1}C$ is non-singular,
  and therefore the matrix $(A - B D^{-1}C) L_{\omega,\phi}^{1/2}$ is itself non-singular. 
  But in the situation of (\ref{newcompactform}), it is well-know that
  the solution $\langle({\underline u}_1,\ldots, {\underline u}_d),\phi\rangle$ is a 
  stochastic process with absolutely continuous law for any test
  function $\phi\not\equiv 0$.
  
  We conclude that
  the law of $\langle u,\phi\rangle$ conditioned to 
  $\langle w,\phi\rangle$, which coincides with the law of    
  $\langle\underline u,\phi\rangle$, is absolutely continuous 
  almost surely with respect to the law of $w$. This is sufficient to conclude that 
  $\langle(u_1,\ldots, u_d, v_1, \ldots, v_q),\phi\rangle$ has
  an absolutely continuous law, which completes the proof.
\qed

\setcounter{equation}{0}  
\section{Example: An electrical circuit}
\label{ElecCirc} 
In this last section we shall present an example of  
linear SDAE's arising from a problem of electrical circuit  
simulation.  

An electrical circuit is a set of devices connected by wires. 
Each device has two or more connection {\em ports.} 
A wire connects two devices at specific ports. Between any 
two ports of a device
there is a {\em flow} (current) and a {\em tension} (voltage drop). 
Flow and tension are supposed to be the same at both ends of a 
wire; thus wires are just physical media for putting together
two ports and they play no other role. 

The circuit topology can be conveniently represented by a network, 
i.e. a set of nodes and a set of directed arcs between nodes, in the 
following way: Each port is a node (taking into account that two
ports connected by a wire collapse to the same node), and any two ports
of a device are joined by an arc. Therefore, flow and tension will be two
quantities circulating through the arcs of the network.
  
It is well known that a network can be univocally described by an  
incidence matrix $A=(a_{ij})$. If we have $n$ nodes and $m$ arcs, 
$A$ is the $m\times n$ matrix defined by
$$
  a_{ij}=
  \left\{\begin{array}{rl}
  +1,\ &\ \text{if arc $j$ has node $i$ as origin}
  \\
  -1,\ &\ \text{if arc $j$ has node $i$ as destiny}
  \\
  0,\ &\ \text{in any other case.}
  \end{array}\right.    
$$

At every node $i$, a quantity $d_i$ (positive, negative or null) of flow 
may be supplied from the outside. This quantity, added to the total
flow through the arcs leaving the node, must equal the total flow   
arriving to the node. This conservation law leads to the system
of equations $Ax=d$, where $x_j$, $j=1,\dots,n$, is the flow through
arc $j$. 

A second conservation law relates to tensions and the cycles formed by
the flows. A {\em cycle} is a set of arcs carrying nonzero flow when 
all external supplies are set to zero. The {\em cycle space} is thus
$\ker A\subset \R^n$. Let $B$ be a matrix whose columns form a basis
of the cycle space, and let $c\in\R^n$ be the vector of externally supplied 
tensions to the cycles of the chosen basis. Then we must impose 
the equalities $B^\top u=c$, where $u_j$, $j=1,\dots,n$, is the tension
through arc $j$. 

Once we have the topology described by a network, we can put into play the
last element of the circuit modelling. Every device has a specific
behaviour, which is described by an equation $\varphi(x,u,\dot x,\dot u)=0$
involving in general flows, tensions, and their derivatives. The system
$\Phi(x,u,\dot x,\dot u)=0$
consisting of all of these equations is called the {\em network characteristic}.
For instance, typical simple two-port (linear) devices are the {\em resistor},
the {\em inductor} and the {\em capacitor}, whose characteristic (noiseless)
equations, which involve only their own arc $j$, 
are $u_j=R x_j$, $u_j=L \dot x_j$, and $x_j=C \dot u_j$,
respectively, for some constants $R,L,C$. 
Also, the {\em current source} ($x_j$ constant) and the {\em voltage source} 
($u_j$ constant) are common devices. 

Solving an electrical circuit therefore means finding the currents $x$ and
voltage drops $u$ determined by the system 
$$
  \left\{\begin{array}{l}
  Ax=d\ 
  \\
  B^\top u=c
  \\
  \Phi(x,u,\dot x,\dot u)=0 
  \end{array}\right.    
$$

\begin{example}  
\label{example}
Let us write down the 
equations corresponding to the circuit called {\em LL-cutset\/} (see \cite{thesisSchein}, pag. 60), 
formed by two inductors and one resistor, which we assume submitted to 
random perturbations, independently for each device. This situation can be modelled,
following the standard procedure described above, by the stochastic system
\begin{equation}
\label{LL-cutset}
\left\{
\begin{array}{l}
x_1=-x_2=x_3 \\
u_1-u_2+u_3=0 \\
u_1=L_1 \dot x_1 + \tau_1 \xi_1 \\
u_2=L_2 \dot x_2 + \tau_2 \xi_2 \\
u_3=R x_3 + \tau_3 \xi_3 
\end{array}
\right.
\end{equation}
where $\xi_1,\xi_2,\xi_3$ are independent white noises,
and $\tau_1,\tau_2,\tau_3$ are non-zero constants which measure the magnitude 
of the perturbations.
With a slight obvious simplification, 
we obtain from
(\ref{LL-cutset})
the following linear SDAE:  
\begin{equation}  
\label{EX2}  
  \left(  
  \begin{matrix}  
  0 & 0 & 0 & 0 \\  
  0 & 0 & 0 & 0 \\  
  0 & 0 & L_1 & 0 \\  
  0 & 0 & 0 & L_2  
  \end{matrix}  
  \right)  
  \left(  
  \begin{matrix}  
  \dot u_1 \\  
  \dot u_2 \\  
  \dot x_1 \\  
  \dot x_2  
  \end{matrix}  
  \right)  
  +
  \left(  
  \begin{matrix}  
  R^{-1} & -R^{-1} & 1 & 0\\  
  -R^{-1} & R^{-1} & 0 & 1 \\  
  -1 & 0 & 0 & 0\\  
  0 & -1 & 0 & 0  
  \end{matrix}  
  \right)  
  \left(  
  \begin{matrix}  
  u_1 \\  
  u_2 \\  
  x_1 \\  
  x_2  
  \end{matrix}  
  \right)  
  =  
  \left(  
  \begin{matrix}  
  0 & 0 & -\tau_3R^{-1} \\  
  0 & 0 & \tau_3R^{-1} \\  
  -\tau_1 & 0 & 0 \\  
  0 & -\tau_2 & 0  
  \end{matrix}  
  \right)  
  \left(  
  \begin{matrix}  
  \xi_1 \\  
  \xi_2 \\  
  \xi_3
  \end{matrix}  
  \right)  
  \ ,  
\end{equation}  
Let us now reduce the equation to KCF. To  
simplify the exposition, we shall fix to 1 the values of $\tau_i$, $R$ and  
$L_i$. (A physically meaningful magnitude for $R$ and $L_i$  
would be of order $10^{-6}$ for the first and of order $10^{4}$ for  
the latter. Nevertheless the structure of the problem does not  
change with different constants.) 
The matrices $P$ and $Q$, providing the desired  
reduction (see Lemma \ref{kronecker}), are 
\[
  P =  
  \left(  
  \begin{matrix}                        \vspace{1mm}
  \frac{1}{2} & -\frac{1}{2} & 1 & -1\\ \vspace{1mm} 
  0 & -1 & 1 & 1 \\                     \vspace{1mm}
  1 & 1 & 0 & 0 \\                      
  -1 & 0 & 0 & 0  
  \end{matrix}  
  \right)  
  \quad,\quad
  Q = \left(  
  \begin{matrix}                                    \vspace{1mm}
  -\frac{1}{4} & -\frac{1}{2} & -\frac{3}{4} & -1\\ \vspace{1mm} 
  \frac{1}{4} & -\frac{1}{2} & -\frac{1}{4} & 0 \\  \vspace{1mm} 
  \frac{1}{2} & 0 & \frac{1}{2} & 0 \\              
  -\frac{1}{2} & 0 & \frac{1}{2} & 0  
  \end{matrix}  
  \right)  
  \quad.
\]
  Indeed, multiplying (\ref{EX2}) by $P$ from the left and defining
  $y=Q^{-1}x$, we arrive to 
\begin{equation}  
\label{EX3}
  \left(  
  \begin{matrix}  
  1 & 0 & 0 & 0 \\  
  0 & 0 & 1 & 0 \\  
  0 & 0 & 0 & 0 \\  
  0 & 0 & 0 & 0  
  \end{matrix}  
  \right)  
  \dot y(t)+  
  \left(  
  \begin{matrix}  
  \frac{1}{2} & 0 & 0 & 0\\  
  0 & 1 & 0 & 0 \\  
  0 & 0 & 1 & 0\\  
  0 & 0 & 0 & 1  
  \end{matrix}  
  \right)  
  y(t)  
  =  
  \left(  
  \begin{matrix}  
  -\tau_1 & \tau_2 & -\tau_3 \\  
  -\tau_1 & -\tau_2 & -\tau_3 \\  
  0 & 0 & 0 \\  
  0 & 0 & \tau_3  
  \end{matrix}  
  \right)  
  \left(  
  \begin{matrix}  
  \xi_1 \\  
  \xi_2 \\  
  \xi_3
  \end{matrix}  
  \right)  
  \ ,  
\end{equation}  
  We see that the matrix $N$ of Section \ref{eu} has here two blocks:
  A single zero in the last position ($\dot y_4$) and a 2-nilpotent block
  affecting $\dot y_2$ and $\dot y_3$.
  We have therefore an index 2 SDAE. 
  From Propositions \ref{AlgPart} and Theorem \ref{Thm43}, 
  we can already say that, when applied to any 
  test function $\phi\neq 0$, 
  the variables $y_4$, $y_2$ and $y_1$, as well as the vectors $(y_1,y_2)$ and $(y_1,y_4)$,  
  will be absolutely
  continuous, whereas $y_3$ degenerates to a point. 
  
  In fact, in this case, we can of course solve completely the system:
  The differential part is the
  one-dimensional classical SDE  
\begin{equation}  
\label{STOex2}  
  \dot y_1  
  +  
  \frac{1}{2} y_1
  =  
  -\tau_1\xi_1+\tau_2\xi_2-\tau_3\xi_3 
  \ ,
\end{equation}  
  and the algebraic part reads simply  
\begin{equation}  
\label{ALex2}  
  \left\{  
  \begin{array}{l}  
  \dot y_3  
  +  
  y_2  
  =  
  -\tau_1\xi_1-\tau_2\xi_2-\tau_3\xi_3 
  \\
  y_3 
  =  
  0  
  \\   
  y_4  
  =  
  \tau_3\xi_3  
  \ .
  \end{array}  
  \right.  
\end{equation}  
  The solution to (\ref{EX3}) can thus be written as
\begin{align*}
y_1(t)  &=
e^{-(t-t_0)/2} 
\Big[y(t_0) + \int_{t_0}^t e^{-(s-t_0)/2} (-\tau_1dW_1+\tau_2dW_2-\tau_3dW_3)(s)\Big]
\\
y_2 &=-\tau_1\xi_1-\tau_2\xi_2-\tau_3\xi_3
\\
y_3 &=0 
\\
y_4 &=\tau_3\xi_3
\ ,
\end{align*}
  where $W_1$, $W_2$, $W_3$ are independent Wiener processes whose generalised 
  derivatives are $\xi_1$, $\xi_2$ and $\xi_3$.
  Multiplying by the matrix $Q$ we finally obtain the value of the original
  variables:
\begin{align*}
x_1(t) &= -x_2(t)  = -\tfrac{1}{2}y_1(t)
\\
u_1 &= -\tfrac{1}{4}y_1-\tfrac{1}{2}y_2-\tfrac{3}{4}y_4
\\
u_2 &= \tfrac{1}{4}y_1-\tfrac{1}{2}y_2 
\ ,
\end{align*}
with $x_1(t_0)=-\frac{1}{4}y_1(t_0)$ a given intensity at time $t_0$.

It is clear that the current intensities, which have
almost surely continuous paths, are much more
regular than the voltage drops, which are only random distributions.

\end{example}


\begin{thebibliography}{10}

\bibitem{MR41:6332}
Donald~A. Dawson.
\newblock Generalized stochastic integrals and equations.
\newblock {\em Trans. Amer. Math. Soc.}, 147:473--506, 1970.

\bibitem{MR0221576}
Xavier Fernique.
\newblock Processus lin\'eaires, processus g\'en\'eralis\'es.
\newblock {\em Ann. Inst. Fourier (Grenoble)}, 17(fasc. 1):1--92, 1967.

\bibitem{MR0435834}
I.~M. Gel{\cprime}fand and N.~Ya. Vilenkin.
\newblock {\em Generalized functions. {V}ol. 4}.
\newblock Academic Press, New York, 1964.

\bibitem{MR90d:34004}
E.~Griepentrog and R.~M{\"a}rz.
\newblock Basic properties of some differential-algebraic equations.
\newblock {\em Z. Anal. Anwendungen}, 8(1):25--41, 1989.

\bibitem{Horn}
Roger~A. Horn and Charles~R. Johnson.
\newblock {\em Matrix Analysis}.
\newblock Cambridge University Press, Cambridge, 1990.

\bibitem{MR1956867}
Jean Jacod and Philip Protter.
\newblock {\em Probability essentials}.
\newblock Universitext. Springer-Verlag, Berlin, second edition, 2003.

\bibitem{MR0461643}
Hui~Hsiung Kuo.
\newblock {\em Gaussian measures in {B}anach spaces}.
\newblock Springer-Verlag, Berlin, 1975.

\bibitem{MR2200233}
David Nualart.
\newblock {\em The {M}alliavin calculus and related topics}.
\newblock Probability and its Applications (New York). Springer-Verlag, Berlin,
  second edition, 2006.

\bibitem{MR97e:34011}
Patrick~J. Rabier and Werner~C. Rheinboldt.
\newblock Classical and generalized solutions of time-dependent linear
  differential-algebraic equations.
\newblock {\em Linear Algebra Appl.}, 245:259--293, 1996.

\bibitem{MR2003f:34002}
Patrick~J. Rabier and Werner~C. Rheinboldt.
\newblock Theoretical and numerical analysis of differential-algebraic
  equations.
\newblock In {\em Handbook of numerical analysis, Vol. VIII}, pages 183--540.
  North-Holland, Amsterdam, 2002.

\bibitem{thesisSchein}
O.~Schein.
\newblock {\em Stochastic differential-algebraic equations in circuit
  simulation}.
\newblock PhD thesis, Technische Universit\"{a}t Darmstadt, 1999.

\bibitem{MR99i:65076}
O.~Schein and G.~Denk.
\newblock Numerical solution of stochastic differential-algebraic equations
  with applications to transient noise simulation of microelectronic circuits.
\newblock {\em J. Comput. Appl. Math.}, 100(1):77--92, 1998.

\bibitem{MR0209834}
Laurent Schwartz.
\newblock {\em Th\'eorie des distributions}.
\newblock Publications de l'Institut de Math\'ematique de l'Universit\'e de
  Strasbourg, No. IX-X. Nouvelle \'edition, enti\'erement corrig\'ee, refondue
  et augment\'ee. Hermann, Paris, 1966.

\bibitem{MR2004f:60133}
Renate Winkler.
\newblock Stochastic differential algebraic equations of index 1 and
  applications in circuit simulation.
\newblock {\em J. Comput. Appl. Math.}, 157(2):477--505, 2003.

\end{thebibliography}
\bibliographystyle{plain}  
  
\def\cprime{$'$}

\end{document}